\documentclass[fleqn,12pt]{article}

\usepackage[headings]{espcrc1}
\usepackage{graphicx}
\usepackage[figuresright]{rotating}
\usepackage{setspace}

\newcommand{\AmS}{{\protect\the\textfont2
  A\kern-.1667em\lower.5ex\hbox{M}\kern-.125emS}}

\newtheorem{theorem}{Theorem}
\newtheorem{lemma}{Lemma}

\newtheorem{conjecture}{Conjecture}

\newtheorem{proposition}{Proposition}

\begin{document}

\title{Tutte's 5-flow conjecture for highly cyclically connected cubic graphs}

\author{Eckhard Steffen \\
	{Paderborn Institute for Advanced Studies in Computer Science and Engineering, 
        Universit\"at Paderborn, Warburger Stra{\ss}e 100, D-33098 Paderborn, Germany} }

\runtitle{5-flows on highly cyclically connected cubic graphs}
\runauthor{E. Steffen}

\maketitle

\begin{abstract}
{In 1954, Tutte conjectured that every bridgeless graph has a nowhere-zero 5-flow. 
Let $\omega(G)$ be the minimum number of odd cycles in a 2-factor of a bridgeless cubic graph $G$. 
Tutte's conjecture is equivalent to its restriction to cubic graphs with $\omega \geq 2$. 
We show that if a cubic graph $G$ has no edge cut with fewer than $ \frac{5}{2} \omega(G) - 3$ 
edges that separates two odd cycles of a minimum 2-factor of $G$, then $G$ has a nowhere-zero 5-flow. 
This implies that if a cubic graph $G$ is cyclically $n$-edge connected and 
$n \geq  \frac{5}{2} \omega(G) - 3$, then $G$ has a nowhere-zero 5-flow.}
\end{abstract}

\section[]{Introduction}

This paper is about flows on finite graphs.
Let $M=(V,E)$ be a graph with vertex set $V$ and 
edge set $E$. Each edge is incident to precisely two different vertices, 
i.e.~multiple edges may occur but there are no loops. 

An {\em orientation} $D$ of $M$ is an assignment of a 
direction to each edge, and for $v \in V,$ 
$D^-(v)$ ($D^+(v)$) is the set of edges whose head 
(tail) is incident to $v$. The oriented graph is denoted by $D(M),$ 
$d_{D(M)}^-(v) = |D^-(v)|$ and $d_{D(M)}^+(v) = |D^+(v)|$ denote the {\em indegree}
and {\em outdegree} of vertex $v$ in $D(M),$ respectively.  

Let $k \geq 2$ be a positive integer and 
$\varphi : E \longrightarrow \{0,1, \dots, k-1\}$ be a function. If for all $v \in V,$

\begin{equation} \label{flow}
\sum_{e \in D^+(v)}\varphi(e) = \sum_{e \in D^-(v)}\varphi(e),
\end{equation} 

then $(D, \varphi)$ is a $k${\em -flow} on $M$. If, in addition, $\varphi(e) \not = 0,$
for all $e \in E,$ then $(D, \varphi)$ is a {\em nowhere-zero} $k$-flow on $M$.
In such a case, we say that $M$ has a nowhere-zero $k$-flow. 

If a graph has a nowhere-zero $k$-flow, then it has one for
every $k' \geq k$. Tutte \cite{Tutte_49} proved 
that a graph $G$ has a nowhere-zero $k$-flow $(D,\varphi)$ if and
only if it has a flow
$(D',\varphi')$ such that for every edge $e$, $|\varphi'(e)|$ is one of 
$1, \dots, k-1$. Thus determining for which number $k$ a graph has a nowhere-zero
$k$-flow is a problem about graphs, not directed graphs.

Tutte \cite{Tutte_54} raised the problem to determine the smallest number $k$
for which a graph has a nowhere-zero $k$-flow, and he formulated 
the 5-Flow Conjecture. 

\begin{conjecture} \label{5fc} \cite{Tutte_54}
Every bridgeless graph has a nowhere-zero 5-flow.
\end{conjecture}

The 5-Flow Conjecture is equivalent to its restriction to cubic graphs, cf.
\cite{Jaeger_88}. By Petersen's theorem, every bridgeless cubic graph $G$ has a 2-factor
and the {\em oddness} $\omega (G)$ is the minimum number of odd cycles in a
2-factor of $G$. Clearly, the oddness must be an even number, and it is well 
known (cf. \cite{Jaeger_88}) that a cubic graph $G$ has a 
nowhere-zero 4-flow if and only if it is edge 3-colorable (i.e. $\omega (G) = 0$). 
Hence the 5-Flow Conjecture is equivalent to its restriction to bridgeless cubic graphs with
$\omega \geq 2$.

Many papers deal with the structure of 
a possible counterexample to the 5-Flow Conjecture. 
A connected graph $G = (V, E)$ that contains two disjoint cycles is 
{\em cyclically $n$-edge connected} if there is no edge cut $E' \subset E$
with fewer than $n$ edges such that two components of $G - E'$ contain cycles.
The maximum number $k$ so that $G$ is cyclically $k$-edge connected
is the {\em cyclic connectivity} of $G$ and it is denoted by $n_G^*$. 
Kochol  \cite{Kochol_04,Kochol_06} showed that the length of a shortest cycle in
a possible minimum counterexample is at least 9, and that it is cyclically 6-edge connected. 
This paper proves the following theorems. 

\begin{theorem} \label{schoenesErgebnis} Every cubic graph $G$ with
cyclic connectivity $n_G^* \geq   \frac{5}{2} \omega(G) - 3$ has a nowhere-zero 5-flow.
\end{theorem}

A {\em minimum} 2-factor of a cubic graph $G=(V,E)$ has precisely $\omega(G)$ odd cycles.
Let $\omega(G) \geq 2$, ${\cal F}_2$ be a minimum 2-factor, and let 
$m_G({\cal F}_2)$ be the maximum number $k$ such 
that there is no edge cut $E' \subset E$ with fewer than $k$ edges such that two
components of $G - E'$ contain odd cycles of ${\cal F}_2$. 
We define 
$m_G^* = \max \{ m_G({\cal F}_2) | {\cal F}_2 \mbox{ is a minimum 2-factor of } G \}$
to be the {\em cyclic factor connectivity} of $G$. For graphs $G$ with $\omega(G) = 0$
define $m_G^* = \infty$. 

Since $n_G^* \leq m_G^*$ Theorem \ref{schoenesErgebnis} is a direct 
consequence of the following theorem.
 
\begin{theorem} \label{Main}
Let $G$ be a bridgeless cubic graph. If $m_G^* \geq \frac{5}{2} \omega(G) - 3$,
then $G$ has a nowhere-zero 5-flow.
\end{theorem}

\section[]{Balanced valuations and flow partitions}

Bondy \cite{Bondy} and Jaeger \cite{Jaeger_75} introduced the concept of  balanced valuations. 
A {\em balanced valuation} of a graph
$M = (V,E)$ is a function $w$ from the vertex set $V$ into the real numbers such that
for all $X \subseteq V$: $| \sum_{v \in X} w(v) | \leq | \partial_M(X) |$, where 
$\partial_M(X)$ is the set of edges with precisely one end in $X$. For $v \in V$
let $d_M(v)$ be the degree of $v$ in the undirected graph $M$.
The following theorem relates integer flows to balanced valuations.

\begin{theorem} \cite{Jaeger_75} \label{Thm_Jaeger_75} 
Let $M=(V,E)$ be a graph with orientation $D$ and $k\geq 3$. Then $M$ 
has a nowhere-zero $k$-flow $(D,\varphi)$ if and only if there 
is a balanced valuation $w$ of $M$ with 
$ w(v) = \frac{k}{k-2}(2d_{D(M)}^+(v) - d_M(v))$, for all $v \in V.$
\end{theorem}

In particular, Theorem \ref{Thm_Jaeger_75} says that a cubic graph $G$ has a nowhere-zero 4-flow 
(nowhere-zero 5-flow) if and only if there is a balanced valuation of $G$ with values in 
$\{ \pm 2\}$ ($\{ \pm \frac{5}{3}\}$).
 
Let $M=(V,E)$ be a multigraph. If $X \subseteq E$, then $M[X]$ denotes the graph whose vertex set 
consists of all vertices of 
edges of $X$ and whose edge set is $X$. Likewise if $X \subseteq V,$ then $M[X]$ is the graph
whose vertex set is $X$ and whose edge set consists of those edges incident to two vertices 
of $X$. In both instances the subgraph $M[X]$ is called the subgraph of $M$ 
{\em induced} by $X$.

Let 
$E_i \subseteq E$, and $(D_i,\varphi_i)$
be flows on $M[E_i]$, $i=1,2$.
The {\em sum} $(D_1,\varphi_1) + (D_2,\varphi_2)$ is the flow $(D, \varphi)$
on $M [E_1 \cup  E_2]$ with orientation

\[ D := D_1|_{ \{e|\varphi_1(e)   \geq  \varphi_2(e)\}}    \cup 
          D_2|_{\{e|\varphi_2(e) > \varphi_1(e)\}} \mbox{, and } \] 

\[ \varphi(e) := \left\{ \begin{array} {l@{\mbox{ }}l}
              \varphi_1(e) + \varphi_2(e) & \mbox{if } e  
              \mbox{ received the same direction
                in } D_1 \mbox{ and } D_2\\
              |\varphi_1(e) - \varphi_2(e)| & \mbox{ otherwise, }
                      \end{array} \right. \]
for $e \in E_1 \cup E_2$.

\vspace{.3cm}
Let $G = (V,E)$ be a bridgeless cubic graph, and ${\cal F}_2$ be a 2-factor of $G$ with 
odd cycles $C_1, C_2, \dots, C_{2 t}$,
and even cycles $C_{2t + 1}, \dots, C_{2t + \ell}$ $(t \geq 0,  \ell \geq 0)$,
and let ${\cal F}_1$ be the complementary 1-factor.

A {\em canonical} 4-coloring of $G$ (with respect to ${\cal F}_2$) colors
the edges of ${\cal F}_1$ with color 1, 
the edges of the even cycles with 2 and 3, alternately, and the
edges of the odd cycles with colors 2 and 3 alternately, except one edge which is colored 0.
Then, there are precisely $2t$ vertices $z_1, z_2, \dots ,z_{2t}$ where color $2$ is missing.

Let $M_G=(V,E(M_G))$ be the graph obtained from $G$ by adding two edges $f_i$ and $f_i'$ 
between $z_{2i-1}$ and $z_{2i}$ for $i = 1, \dots, t$. Extend $c$ to a proper edge coloring of 
$M_G$ by coloring 
$f'_i$ with color 2 and $f_i$ with color 4. Let $C'_1, \dots C'_s$ be the 2-factor of $M_G$
induced by the edges of colors $1$ and $2$ ($s \geq 1$), and for
$i= 1, \dots, t$  let $C''_i$ be the 2-cycle induced by the edges $f_i$ and $f_i'$. 
We construct a nowhere-zero 4-flow on $M_G$ as follows:

For $1 \leq i \leq 2 t + \ell$ let $(D_i, \varphi_i)$ be a nowhere-zero flow on the directed cycle 
$C_i$ with $\varphi_i(e) = 2$ for all $e \in E (C_i)$.

For  $1 \leq i \leq s$ let $(D'_i, \varphi'_i)$ be a nowhere-zero flow on the directed cycle $C'_i$ with 
$\varphi_i' (e) = 1$ for all $e \in E (C'_i)$.

For $1 \leq i \leq t$ let $(D''_i, \varphi''_i)$ be a nowhere-zero flow on the directed cycle $C''_i$ 
(choose $D''_i$ such that $f'_i$ receives the same direction as in $D'_i$) with $\varphi''_i (e) = 1$ 
for all $e \in \{f_i, f'_i \}$.
Then

\begin{eqnarray} \label{1.1}
(D, \varphi) = \sum_{i=1}^{2 t + \ell} (D_i, \varphi_i) 
						+ \sum_{i=1}^s (D'_i, \varphi'_i) 
							+ \sum^t_{i=1} (D''_i, \varphi''_i)
\end{eqnarray}

is a nowhere-zero 4-flow on $M_G$.

By Theorem \ref{Thm_Jaeger_75}, there is a balanced valuation $w' (v) = 2 (2 d^+_{D(M_G)}(v) - d_{M_G}(v))$
of $M_G$. It holds that $|2 d^+_{D(M_G)}(v) - d_{M_G}(v)| = 1$, and hence
$w'(v) \in \{ \pm 2 \}$ for all $v \in V$. 
The vertices of $M_G$ (and therefore of $G$ as well) are
partitioned into two classes $A = \{ v | w'(v) = -2\}$ and $B = \{v | w'(v) = 2\}$. 
Call the elements of $A$ ($B$) the white (black) vertices of $M_G$ and of $G$, respectively. 

Let $G=(V,E)$ be a bridgeless cubic graph. A partition of $V$ into two classes 
$A$ and $B$ constructed as above, and using 
a 2-factor ${\cal F}_2$, a canonical 4-coloring $c$ of $G$, the 
4-flow $(D,\varphi)$ on $M_G$ and the induced balanced valuation $w'$ of $M_G$   
is called a {\em flow partition of $G$}, and it is 
denoted by $P_G(A,B) = P_G(A,B,{\cal F}_2,c,(D,\varphi),w')$. 
If we refer to a special 
2-factor ${\cal F}_2$, we say $P_G(A,B)$ is a flow partition of $G$ with respect to
${\cal F}_2$.
For $X \subseteq V$ let $A_X = A \cap X$ ($B_X = B \cap X$) be the set of the white (black) 
vertices of $X$, and $a_X = |A_X|$, $b_X = |B_X|$. If we consider the vertex set $V(F)$
of a subgraph $F$ of a graph $G$ we also write
$a_F$ instead of $a_{V(F)}$ ($b_F$ instead of $b_{V(F)}$).

We will prove some properties of flow partitions of cubic graphs.
The following Lemma is a direct consequence of the construction of $(D,\varphi)$ on
$M_G$. 

\begin{lemma} \label{verschiedene2}
Let $P_G(A,B,{\cal F}_2,c,(D,\varphi),w')$ be a flow partition of a 
bridgeless cubic graph $G=(V,E)$, and $xy = e \in E$. If the 
canonical 4-coloring $c$ colors $e$ with $1$ or $2$, then
$x$ and $y$ belong to different classes, i.e. $x \in A$ if and only if $y \in B$.
\end{lemma}

\begin{lemma} \label{Schranke} 
Let $G = (V,E)$ be a cubic bridgeless graph and $P_G(A,B)$ be a flow partition
with respect to a 2-factor ${\cal F}_2$. Let 
$S \subseteq V$ be a set of vertices such that the induced subgraph
$G[S]$ is connected, $n$ be the
number of edges which have to be removed from $G[S]$ to obtain a spanning tree
of $G[S]$, and let $n_o$ be the number of odd cycles of ${\cal F}_2$ which are 
subgraphs of $G[S]$. Then $b_S \leq 4a_S + 3 - 3n + n_o$.
\end{lemma}

Proof.  
Let $P_G(A,B) = P_G(A,B,{\cal F}_2,c,(D,\varphi),w')$ and 
$F$ be a connected subgraph of ${\cal F}_2$. We show:

1) If $F$ is an even cycle, then $b_F = a_F$.

2) If $F$ is an odd cycle, then $b_F \leq  a_F + 1$.

3) If $F$ is a path, then $b_F \leq  a_F + 3$.

Items 1) and 2) follow from Lemma \ref{verschiedene2} directly.
We distinguish two cases to prove 3).

Case 1: The edges of $F$ are colored with colors $2$ and $3$. 

If $|E(F)| = 2l+1$, then at least $l$ edges 
are colored with color $2$. Thus Lemma \ref{verschiedene2} implies that $a_F \geq l$. Since $|V(F)| = 2l+2$ and 
$b_F = 2l+2 - a_F$ it follows that $b_F \leq a_F +2$.

If $|E(F)| = 2l$, then $l$ edges are colored with color $2$. Thus Lemma \ref{verschiedene2}  implies that $a_F = l$. 
Since $|V(F)| = 2l+1$ it follows that $b_F = a_F +1$.

Case 2: $F$ contains an edge of color $0$. 

By the definition of the coloring there is precisely one edge of color $0$.

If the length of $F$ is odd, say $2l+1$, the first and the last edge of $F$
are colored differently, and there are $l$ edges of color $2$. Thus 
Lemma \ref{verschiedene2} implies that $a_F \geq l$. Since $|V(F)| = 2l+2$ it follows that 
$b_F \leq a_F + 2$.

If $|E(F)| = 2l$, then at least $l-1$ edges are colored $2$. Thus Lemma \ref{verschiedene2}  implies 
that $a_F \geq l-1$. Since $|V(F)| = 2l+1$ it follows that $b_F \leq a_F + 3$. $\circ$

Let $E_1$ be the set of edges of color 1 of $G[S]$. By Lemma \ref{verschiedene2},
$|E_1| \leq a_S$. Let $E_1^- \subset E_1$ be a set of edges so that $G[S] - E_1^-$ is
connected and no edge of color 1 (in $G[S] - E_1^-$) is contained in a cycle. Each cycle of $G[S] - E_1^-$
is a cycle of ${\cal F}_2$. Remove from each cycle precisely one edge of color 2 to 
obtain a spanning tree of $G[S]$. Let $E_2^-$ be the set of these removed edges of color 2. 
With $n_i = |E_i^-|$ ($i = 1,2$) it follows that $n = n_1 + n_2$. 

Let  $Z_0, \dots , Z_{a'}$ be the components of $G[S] - E_1$, and
$a_i$ ($b_i$) be the number of white (black) vertices in $Z_i$, $i = 0, \dots, a'$.
Each component is either a cycle of ${\cal F}_2$ or a subpath of a cycle of ${\cal{F}}_2$.
The number of components is smaller than or equal to 1 plus the number of edges of color 1
in $G[S] - E_1^-$, therefore $a' \leq a_S - n_1$. Furthermore $\sum_{i=0}^{a'}a_i = a_S$.

For $i \in I_P = \{0, 1, \dots ,a' - n_2\}$ let $Z_i$ be a path,
for $i \in I_C^o = \{ a' - n_2 +1, \dots ,a' - n_2 + n_o\}$ let $Z_i$ be an odd cycle,
and for $i \in I_C^e = \{a' - n_2 + n_o + 1, \dots , a'\}$ let $Z_i$ be an even cycle.
Then it follows with $a' \leq a_S - n_1$ that

{\onehalfspacing
$b_S 	= 	\sum_{i \in I_P} b_i	+ \sum_{i \in I_C^o} b_i +  \sum_{i \in I_C^e} b_i \\
	\leq 	\sum_{i \in I_P} (a_i+3) + \sum_{i \in I_C^o} (a_i + 1) + \sum_{i \in I_C^e} a_i \\
	=	3(a'-n_2+1) + n_o +  \sum_{i=0}^{a'} a_i \\
	\leq  3(a_S -(n_1+n_2) + 1) + n_o + a_S \\
	= 4a_S + 3 -3n + n_o$. \hfill\ {\raisebox{0.8ex}{\framebox{}}}\par\bigskip}%

We finish this section with the following lemma.

\begin{lemma} \label{hinreichendeBedingung} Let $P_G(A,B)$
be a flow partition of a cubic bridgeless graph $G = (V,E)$. 
Let $S \subseteq V$ be a set of vertices such that the induced subgraph
$G[S]$ is connected, and $n$ be the
number of edges which have to be removed from $G[S]$ to obtain a spanning tree $T$
of $G[S]$. Assume $a_S \leq b_S$, then
$b_S \leq 4a_S + 3 - 3n$ if and only if  $\frac{5}{3}(b_S - a_S) \leq |\partial_G(S)|$.
\end{lemma}

Proof. Consider a spanning tree $T = (S,E(T))$ of $G[S]$ and let
$T_i = \{v | v \in S$ and  $d_{T} (v) = i \}$, for $i = 1, 2$.
Then $|\partial_G(S)| + 2n = 2 |T_1| + |T_2|$ and \vspace{.2cm} 

$|S|-1 = |E(T)| = \frac{1}{2} (3 (|S| - (|T_1| + |T_2|)) + 2 |T_2| + |T_1|) 
=  \frac{1}{2} (3 |S| - |\partial_G(S)| - 2n).$ \vspace{.2cm} 

Since $|S| = a_S + b_S$ it
follows that $|\partial_G(S)| = a_S + b_S + 2-2n$, and hence
$\frac{5}{3}(b_S - a_S) \leq |\partial_G(S)|$ is equivalent 
to $b_S \leq 4a_S + 3 - 3n$.  \hfill\ {\raisebox{0.8ex}{\framebox{}}}\par\bigskip

\section[]{Proof of Theorem \ref{Main}}

Let $G = (V,E)$ be a bridgeless cubic graph with oddness $\omega$. 
If $\omega \in \{0,2\}$, then $G$ has a nowhere-zero 5-flow, cf. \cite{Jaeger_88}. Thus we
may assume that $\omega \geq 4$. 

Let ${\cal F}_2$ be a minimum 2-factor of $G$ with 
$m_G({\cal F}_2) = m_G^* \geq \frac{5}{2} \omega - 3$.
Let $P_G(A,B)= P_G(A,B,{\cal F}_2,c,(D,\varphi),w')$ be a flow partition of $G$
with respect to ${\cal F}_2$. Let
$w : V \rightarrow \{\pm \frac{5}{3}\}$ be a function with 
$w(v) = - \frac{5}{3}$ if $v \in A$ and $w(v) = \frac{5}{3}$ if $v \in B$.
We will show that $w$ is a balanced valuation of $G$. Then it follows from Theorem \ref{Thm_Jaeger_75} 
that $G$ has a nowhere-zero 5-flow.

Assume to the contrary that $w$ is not balanced. Then there is $S \subseteq V$ with 
\begin{eqnarray} \label{widerspruch}
| \sum_{v \in S} w(v)| > |\partial_G (S)|.
\end{eqnarray} 

If $S=V$, then $| \sum_{v \in S} w (v)| = 0 = |\partial_G (S)|$, and
therefore $S$ is a proper subset of $V.$ 
Let $S$ be of minimum order, so we may assume that $G[S]$ is connected, and 
without loss of generality $b_S \geq a_S$. 
With $k = b_S-a_S$ equation (\ref{widerspruch}) becomes
\begin{eqnarray} \label{bv}
\frac{5}{3}k > |\partial_G(S)|.  
\end{eqnarray}

We show

\begin{proposition} \label{cut}
$ |\partial_G (S)| \leq \frac{5}{2} \omega - 4$; in particular 

{\onehalfspacing
$1)$ $ |\partial_G (S)| \leq \frac{5}{2} \omega - 4$, if $|\partial_G (S)| \equiv 1 \mbox{ mod } 5$,  

$2)$ $ |\partial_G (S)| \leq \frac{5}{2} \omega - 8$, if $|\partial_G (S)| \equiv 2 \mbox{ mod } 5$, 

$3)$ $ |\partial_G (S)| \leq \frac{5}{2} \omega - 7$, if $|\partial_G (S)| \equiv 3 \mbox{ mod } 5$, 

$4)$ $ |\partial_G (S)| \leq \frac{5}{2} \omega - 11$, if $|\partial_G (S)| \equiv 4 \mbox{ mod } 5$,  

$5)$ $ |\partial_G (S)| \leq \frac{5}{2} \omega - 15$, if $|\partial_G (S)| \equiv 0 \mbox{ mod } 5$. 

}
\end{proposition}

Proof. For $i = 0,1,2,3$ let $E_i \subset E$ be the set of the edges of color $i$ in $G$
and let $c_i = |\partial_G (S) \cap  E_i|$. The edges of color 1 form a 
1-factor of $G$. Thus Lemma \ref{verschiedene2} implies
that $k = c_1$ and hence $c_1 > \frac{3}{5}|\partial_G(S)|$ by equation (\ref{bv}). 

Let $l_S^a$ ($l_S^b$) be the number of white (black) vertices of $S$ where color 2 is missing,
with respect to $c$. 
Let $l = |l_S^b - l_S^a|$. From $0 \leq l_S^a$, $l_S^b \leq \frac{1}{2} \omega$ it follows that $l \leq \frac{1}{2} \omega $, 
and Lemma \ref{verschiedene2} implies that $k \leq c_2 + l$. Hence
$c_2 + \frac{1}{2} \omega \geq k > \frac{3}{5}|\partial_G(S)|$. 

1) If $|\partial_G(S)| \equiv 1 \mbox{ mod } 5$, say $|\partial_G(S)| = 5m+1$, then
it follows that $c_1 \geq 3m+1$ and therefore $c_2 \leq 2m$. 
Thus $\frac{1}{2} \omega \geq 3m+1 - c_2 \geq 3m+1 - 2m = m+1$ and hence $\frac{5}{2}\omega - 4 \geq |\partial_G(S)|$. 

2) can be proved analogously.

3) If $|\partial_G(S)| \equiv 3 \mbox{ mod } 5$, say $|\partial_G(S)| = 5m+3$, then it follows that $c_1 \geq 3m + 2$ 
and therefore $c_2 \leq 2m + 1$. 

If $c_2 = 2m+1$, then $c_1 \leq |\partial_G(S)| - c_2 = 3m+2$ and hence $c_1 = 3m+2$ and $c_0 = c_3 = 0$. 
Let $X$ be the set of vertices of $G[S]$ which are incident (in $G$) to an edge of $|\partial_G (S) \cap  E_2|$, 
and $Y$ be the set of vertices which are incident to an edge 
of color $0$ in $G[S]$.
Color 2 or 3 is missing on each vertex of $X \cup Y$ and $Z = X \cap Y$ consists of those vertices 
of $G[S]$ where both colors, 2 and 3, are missing. Each vertex $z$ of 
$Z$ is incident to an edge $e=zz'$ of color $0$ in $G[S]$. 
Furthermore, color $2$ is missing and color 3 appears at $z'$. 
Therefore, for each vertex of $z \in Z$ there is precisely one vertex $z'$ in $G[S]$ where only
color 2 is missing. Since $|X|=c_2$ is odd and $c_0 = c_3 = 0$ it follows that 
the total number of vertices of $G[S]$ where either color 2 or color 3 is missing is odd.
This is a contradiction, since every path induced by edges of colors 2 and 3 in $G[S]$ 
has precisely two end vertices in $G[S]$. 

Therefore $c_2 \leq 2m$ and hence  
$c_2 + \frac{1}{2} \omega \geq 3m + 2$ implies that $\frac{1}{2} \omega \geq 3m+2 - 2m = m+2$.
Thus $\frac{5}{2} \omega - 7 \geq 5m + 3 = |\partial_G(S)|$.

Items 4) and 5) can be proved  analogously to 3). $\circ$

Since $G$ has no edge cut with fewer than $\frac{5}{2} \omega -3$ edges
that separates two odd cycles of ${\cal F}_2$ it follows 
with Proposition \ref{cut} that $n_o = 0$. Hence $b_S \leq 4a_S + 3 -3n$ by Lemma \ref{Schranke} and therefore
$\frac{5}{3}k \leq |\partial_G(S)|$ by Lemma \ref{hinreichendeBedingung}. This contradicts equation (\ref{bv})
and completes the proof. \hfill {\raisebox{0.8ex}{\framebox{}}}\par\bigskip

\section[]{Remarks on $r$-flows}

The notion of nowhere-zero flows can be extended to rational numbers. 
Let $1 \leq p \leq q$ be integers, and let 
$\varphi$ be a function from the edge set $E$ of the directed graph $G=(V,E)$ (with orientation $D$)
into the rational numbers.
$(D,\varphi)$ is a {\em nowhere-zero $\frac{q}{p} + 1$-flow} on $G=(V,E)$ if 
$1 \leq \varphi(e) \leq \frac{q}{p}$ for all $e \in E$ and equation (\ref{flow}) is
satisfied for all $v \in V.$ The {\em circular flow number} $F_c(G)$ of $G$ is
the minimum number $r$ such that $G$ has a nowhere-zero $r$-flow. 

Seymour \cite{Seymour_81} proved that every bridgeless graph has a
nowhere-zero 6-flow. Some methods of this paper can be extended to the study of 
nowhere-zero $r$-flows on graphs. For instance, it can be proved that $F_c(G) < 6$ for all bridgeless 
cubic graphs $G$ with $m_G^* \geq \frac{3}{2} \omega(G) + 1$.

\end{document}